\documentclass[12pt,fleqn]{amsart}

\usepackage[utf8]{inputenc}
\usepackage[T1]{fontenc}

\usepackage{a4wide}
\usepackage{scrextend}

\usepackage{amssymb,amsthm,latexsym}
\usepackage{mathrsfs}
\usepackage{dsfont}
\usepackage[mathscr]{euscript}
\usepackage{esint}
\usepackage{newtxtext,newtxmath}
\usepackage{bbm}
\usepackage{enumitem}

\usepackage[usenames,dvipsnames,x11names]{xcolor}
\usepackage{url}
\usepackage{microtype}
\usepackage{hyperref}
\usepackage{cleveref}
\hypersetup{colorlinks=true,linkcolor=blue,citecolor=blue,urlcolor=blue}

\renewcommand{\le}{\leqslant}
\renewcommand{\ge}{\geqslant}

\newcommand{\HS}{\mathrm{HS}}
\newcommand{\Law}{\mathcal{L}}  
\newcommand{\cL}{\mathcal{L}}   
\newcommand{\cH}{\mathcal{H}}
\newcommand{\cR}{\mathcal{R}}
\DeclareMathOperator{\Dom}{Dom}
\DeclareMathOperator{\ran}{Ran}

\newcommand{\R}{\mathbb{R}}
\newcommand{\N}{\mathbb{N}}
\newcommand{\E}{\mathbb{E}}
\renewcommand{\P}{\mathbb{P}}

\newtheorem{theorem}{Theorem}[section]
\newtheorem{lemma}[theorem]{Lemma}
\newtheorem{corollary}[theorem]{Corollary}
\newtheorem{assumption}[theorem]{Assumption}
\newtheorem{proposition}[theorem]{Proposition}

\newtheorem{example}[theorem]{Example}

\newtheorem*{theoremA}{Theorem A}
\newtheorem*{theoremB}{Theorem B}
\newtheorem*{theoremC}{Theorem C}

\theoremstyle{definition}
\newtheorem*{remark*}{Remark}
\newtheorem{remark}[theorem]{Remark}

\numberwithin{equation}{section}

\title[Law equivalence for Ornstein--Uhlenbeck dynamics]{Law equivalence for Ornstein--Uhlenbeck dynamics\\
driven by L\'evy noise}
\author[T.~Kania]{Tomasz Kania}
\address[T.~Kania]{Mathematical Institute\\Czech Academy of Sciences\\\v Zitn\'a 25 \\115 67 Praha 1\\Czech Republic  and  Institute of Mathematics and Computer Science\\ Jagiellonian University\\ {\L}ojasiewicza 6, 30-348 Krak\'{o}w, Poland
}
\email{kania@math.cas.cz, tomasz.marcin.kania@gmail.com}
\date{\today}
\thanks{RVO: 67985840.}

\subjclass[2020]{60H15, 60G51, 47D06, 60J25}
\keywords{Ornstein--Uhlenbeck process, L\'evy process, Cameron--Martin space, Gaussian measures, law equivalence, analytic semigroups}

\begin{document}

\begin{abstract}
For stochastic partial differential equations driven by L\'evy noise, understanding when changes in the drift operator preserve the law of the solution is fundamental to filtering, control, and simulation. We extend law-comparison results for Ornstein--Uhlenbeck processes from bounded drift operators to generators of $C_0$-semigroups on a separable Hilbert space. The argument separates the problem into a Gaussian channel and a jump--drift channel. The Gaussian channel is governed by an inverse-covariance Hilbert--Schmidt perturbation condition. The jump--drift channel is handled by a directional Cameron--Martin condition, formulated conditionally on the jump path; no unconditional Novikov estimate is needed for this step. We prove that these hypotheses give absolute continuity of path laws on the Skorohod space, and equivalence whenever the Gaussian channels are equivalent. For purely jump noise we prove a rigidity phenomenon: absolute continuity forces the two solutions to coincide. For analytic semigroups and compound Poisson jumps we give verifiable sufficient conditions in terms of fractional smoothing and the boundedness of $Q^{-1/2}(\widetilde A-A)A^{-\beta}$. Diagonal examples show the sharp role of the Cameron--Martin and inverse-covariance requirements.
\end{abstract}
\maketitle

\section{Introduction}

The equivalence of path laws for solutions to stochastic differential equations is a classical theme in stochastic analysis, with applications ranging from filtering and statistical inference to the numerical simulation of stochastic partial differential equations (SPDEs). In finite dimensions, Girsanov's theorem gives a robust method for comparing diffusions with different drifts. In infinite-dimensional Hilbert spaces the situation is subtler: Gaussian measures are quasi-invariant only along their Cameron--Martin spaces, and jump noise cannot usually be smoothed away by a change of Gaussian measure.

Let $H$ be a separable Hilbert space and let $L$ be an $H$-valued L\'evy process with L\'evy--It\^o decomposition
\[
        L_t=bt+W_t+Z_t,
\]
where $b\in H$, $W$ is a $Q$-Wiener process and $Z$ is the jump part, independent of $W$. If $A$ generates a $C_0$-semigroup $S=(S(t))_{t\ge0}$, the corresponding Ornstein--Uhlenbeck process is the mild solution
\begin{equation}\label{eq:OU_equation}
        X^A(t)=\int_0^t S(t-s)\,dL_s,
        \qquad 0\le t\le T.
\end{equation}
Given another generator $\widetilde A$ with semigroup $\widetilde S$, we ask when the laws of $X^A$ and $X^{\widetilde A}$ on the Skorohod space $D([0,T];H)$ are absolutely continuous or equivalent.

For bounded drift operators, Bartosz and Kania \cite{BartoszKania2019} obtained law equivalence in the presence of a non-degenerate Wiener component, together with a rigidity theorem in the purely jump case. The present paper treats the generator setting. The principal correction one must make in passing to unbounded generators is that all Gaussian and Cameron--Martin conditions must be measured in the inverse-covariance direction. Thus strict positivity of the eigenvalues of $Q$ means that $Q$ is injective; it does not make $Q^{-1/2}$ bounded or everywhere defined.

Our main results are as follows. Precise statements are given in \Cref{sec:statements,sec:analytic}.

\begin{theoremA}[Mixed Gaussian--jump comparison]
Let $A$ and $\widetilde A$ generate $C_0$-semigroups. Suppose that the Gaussian Ornstein--Uhlenbeck law generated by $A$ is absolutely continuous with respect to the one generated by $\widetilde A$, and that the jump--drift discrepancy
\[
        \int_0^t(S-\widetilde S)(t-s)\,dZ_s
        +\int_0^t(S-\widetilde S)(t-s)b\,ds
\]
belongs, conditionally on the jump path, to the Cameron--Martin space of the $\widetilde A$-Gaussian channel. Then
\[
        \Law(X^A)\ll \Law(X^{\widetilde A})
        \quad\hbox{on }D([0,T];H).
\]
If the two Gaussian Ornstein--Uhlenbeck laws are equivalent, then the two full L\'evy-driven laws are equivalent.
\end{theoremA}

\begin{theoremB}[Pure-jump rigidity]
If $L=Z$ is purely discontinuous with no deterministic drift and no Gaussian part, then
\[
        \Law(X^A)\ll \Law(X^{\widetilde A})
        \quad\Longrightarrow\quad
        X^A=X^{\widetilde A}\quad\hbox{a.s. on }[0,T].
\]
\end{theoremB}

\begin{theoremC}[Analytic sufficient conditions]
Assume that $A$ and $\widetilde A$ generate bounded analytic semigroups and put $K=\widetilde A-A$. If, for some $\beta<\tfrac12$, the operator $Q^{-1/2}KA^{-\beta}$ extends boundedly to $H$, then the Gaussian comparison condition holds. If, in addition, the jump part is compound Poisson, the directional Cameron--Martin condition is satisfied by an explicit representative. Hence the conclusions of Theorem~A apply.
\end{theoremC}

The proof of Theorem~A is a conditional Cameron--Martin argument. For a fixed jump path $z$, the two conditional laws are translates of the Gaussian laws by deterministic paths. The translate difference is precisely the jump--drift discrepancy. This point is important: one must first compare the two Gaussian channels and only then translate by the discrepancy. Doing so avoids an erroneous identification of the conditional Gaussian part.

The paper is organised as follows. \Cref{sec:preliminaries} recalls Skorohod-space measurability, Cameron--Martin shifts for Gaussian convolutions, Duhamel's formula and the Gaussian comparison criterion. \Cref{sec:statements} contains the main abstract hypotheses and theorems. \Cref{sec:analytic} gives analytic-semigroup and compound-Poisson sufficient conditions. \Cref{sec:CMfail} presents diagonal examples. \Cref{sec:proofs} contains the proofs of the main theorems. We conclude with open problems.

\section{Preliminaries}\label{sec:preliminaries}

\subsection{Skorohod space and measurable jump reconstruction}
We work on $D_{H,T}:=D([0,T];H)$, the Polish space of c\`adl\`ag $H$-valued functions on $[0,T]$ with the Skorohod $J_1$ topology. Its Borel $\sigma$-algebra is generated by point evaluations $f\mapsto f(t)$. We use the following measurable reconstruction facts from \cite{BartoszKania2019}.

\begin{proposition}[Measurable jump reconstruction]\label{prop:jump_reconstruction}
There exist Borel maps $\mathcal Z:D_{H,T}\to D_{H,T}$ and, for each Borel set $E\subset H\setminus\{0\}$ bounded away from $0$, maps $Z^1_E:D_{H,T}\to D_{H,T}$ such that:
\begin{enumerate}[label=(\roman*)]
\item if $f$ has the law of a L\'evy process, then $\mathcal Z(f)$ equals its pure jump part almost surely;
\item $Z^1_E(f,t)=\sum_{s\le t,\,\Delta f(s)\in E}\Delta f(s)$ for $t\in[0,T]$;
\item for a generator $B$ with semigroup $S_B$, the map
\[
        f\mapsto \int_0^\cdot S_B(\cdot-s)\,d\mathcal Z(f,s)
\]
is Borel from $D_{H,T}$ to $D_{H,T}$.
\end{enumerate}
\end{proposition}

We shall also use the elementary product-measure lemma below.

\begin{lemma}[{\cite[Lem.~2.2]{BartoszKania2019}}]\label{lem:product_ac}
Let $(X_1,Y)$ and $(X_2,Y)$ be pairs of random variables in a product space with $X_1,X_2$ independent of $Y$. If $\Law(X_1)\ll \Law(X_2)$, then $\Law(X_1,Y)\ll \Law(X_2,Y)$. If $\Law(X_1)\sim\Law(X_2)$, then $\Law(X_1,Y)\sim\Law(X_2,Y)$.
\end{lemma}

\subsection{The inverse covariance and Cameron--Martin shifts}\label{subsec:CM}
Throughout the Gaussian part of the paper $Q$ is a non-negative trace-class operator on $H$. When its eigenvalues are strictly positive, we write
\[
        Qe_n=q_ne_n,\qquad q_n>0,
\]
with respect to an orthonormal basis $(e_n)$. The operator $Q^{1/2}$ is bounded and Hilbert--Schmidt, whereas
\[
        Q^{-1/2}x=\sum_{n\ge1}q_n^{-1/2}\langle x,e_n\rangle e_n
\]
is an unbounded densely defined operator with domain
\begin{equation}\label{eq:domain-qminus}
        \Dom(Q^{-1/2})=
        \left\{x\in H:\sum_{n\ge1}\frac{|\langle x,e_n\rangle|^2}{q_n}<\infty\right\}.
\end{equation}
This is the Cameron--Martin space of the one-time Gaussian law $N(0,Q)$, with norm $\|x\|_Q=\|Q^{-1/2}x\|$.

Let $B$ generate a $C_0$-semigroup $S_B$. The Cameron--Martin operator associated with the Gaussian convolution
\[
        Y_B(W)(t):=\int_0^t S_B(t-s)\,dW_s
\]
is
\begin{equation}\label{eq:KBdef}
        \cR_{B,T}u(t):=\int_0^t S_B(t-s)Q^{1/2}u(s)\,ds,
        \qquad u\in L^2(0,T;H).
\end{equation}
We write
\[
        \cH_{B,T}:=\ran \cR_{B,T}
\]
with the usual Cameron--Martin norm, namely the quotient norm inherited from $L^2(0,T;H)$.

\begin{theorem}[Cameron--Martin theorem for OU convolutions]\label{thm:CM_OU}
Let $B$ generate a $C_0$-semigroup and let $h\in\cH_{B,T}$. Then
\[
        \Law\big(Y_B(W)+h\big)\sim \Law\big(Y_B(W)\big)
        \quad\hbox{on }D_{H,T}.
\]
Equivalently, if $h=\cR_{B,T}u$ for some $u\in L^2(0,T;H)$, the shift is implemented by the usual Cameron--Martin density. If $u$ is random but independent of $W$, the same conclusion holds after conditioning on $u$; no unconditional Novikov estimate is needed in this conditional situation.
\end{theorem}

\begin{remark}\label{rem:no-novikov}
A Novikov condition is required for general adapted shifts depending on the driving Wiener process. In this paper the jump--drift shift is $\sigma(Z)$-measurable and independent of $W$, so the conditional Cameron--Martin theorem is the natural tool. This is also the reason why the covariance must enter through $Q^{-1/2}$ rather than through $Q^{1/2}$.
\end{remark}

\subsection{Duhamel formula for fractionally bounded perturbations}
Let $A$ and $\widetilde A$ generate $C_0$-semigroups $S$ and $\widetilde S$, and set $K:=\widetilde A-A$.

\begin{theorem}[Duhamel formula {\cite[Cor.~III.1.7, Thm.~III.3.14]{EngelNagel2000}}]\label{thm:Duhamel}
If $K\in\cL(H)$, then for $t\ge0$,
\begin{equation}\label{eq:Duhamel}
        \widetilde S(t)-S(t)=\int_0^t \widetilde S(t-s)K S(s)\,ds
\end{equation}
in the uniform operator topology. If $A$ generates a bounded analytic semigroup and there exists $\beta\in(0,1)$ with $KA^{-\beta}\in\cL(H)$, then \eqref{eq:Duhamel} remains valid in the strong operator topology.
\end{theorem}

For bounded analytic semigroups one has the smoothing estimate
\begin{equation}\label{eq:analytic-smoothing}
        \|A^\beta S(t)\|\le C_\beta t^{-\beta},
        \qquad 0<t\le1,
        \quad \beta\in(0,1),
\end{equation}
see \cite[Prop.~II.4.6]{EngelNagel2000}. As usual, a harmless spectral shift may be made before taking fractional powers.

\subsection{Gaussian OU comparison under inverse-covariance Hilbert--Schmidt perturbations}
The following form of the Gaussian comparison criterion is the one used below; it is a standard consequence of the Hilbert-space Girsanov theorem and the Gaussian comparison results of Peszat and Brze\'zniak--van Neerven.

\begin{theorem}[Gaussian OU comparison {\cite{Loges1984,Peszat1992,BrzVanN2000}}]\label{thm:GaussianEquiv}
Let $A$ generate a $C_0$-semigroup $S$, and let $\widetilde A=A+K$ generate $\widetilde S$. Assume that, for each $t>0$, the operator $KS(t)Q^{1/2}$ takes values in $\Dom(Q^{-1/2})$ and that, for every $T>0$,
\begin{equation}\label{eq:Gaussian-HS}
        \int_0^T\big\|Q^{-1/2}KS(t)Q^{1/2}\big\|_{\HS}^2\,dt<\infty .
\end{equation}
Then the Gaussian Ornstein--Uhlenbeck laws satisfy
\[
        \Law(Y_A(W))\sim\Law(Y_{\widetilde A}(W))
        \quad\hbox{on }D_{H,T}.
\]
\end{theorem}

\begin{remark}\label{rem:inverse-covariance}
Condition \eqref{eq:Gaussian-HS} is an inverse-covariance condition. The weaker-looking expression $\|KS(t)Q^{1/2}\|_{\HS}$ controls $KX$ in the ambient Hilbert norm, but not in the Cameron--Martin norm of the noise. In infinite dimension this distinction is essential because $Q^{-1/2}$ is generally unbounded.
\end{remark}

\section{Main hypotheses and theorems}\label{sec:statements}

Write
\[
        X^B=Y_B(W)+J_B(Z)+B_B
\]
where, for a generator $B$ with semigroup $S_B$,
\begin{align*}
Y_B(W)(t)&:=\int_0^t S_B(t-s)\,dW_s, &\text{(Gaussian part)},\\
J_B(Z)(t)&:=\int_0^t S_B(t-s)\,dZ_s, &\text{(jump part)},\\
B_B(t)&:=\int_0^t S_B(t-s)b\,ds. &\text{(deterministic drift)}
\end{align*}

\begin{assumption}[Directional Cameron--Martin condition, $A\to\widetilde A$]\label{ass:CM}
Define the jump--drift discrepancy
\begin{equation}\label{eq:Delta-Jb}
\Delta^{Jb}_{A\to\widetilde A}(t)
:=\int_0^t (S-\widetilde S)(t-s)\,dZ_s
  +\int_0^t (S-\widetilde S)(t-s)b\,ds,
\qquad 0\le t\le T.
\end{equation}
There exists a $\sigma(Z)$-measurable random variable
\[
        U^{A\to\widetilde A}:\Omega\to L^2(0,T;H),
\]
independent of $W$, such that
\begin{equation}\label{eq:CM_representation}
        \Delta^{Jb}_{A\to\widetilde A}(t)
        =\int_0^t \widetilde S(t-s)Q^{1/2}U^{A\to\widetilde A}_s\,ds,
        \qquad 0\le t\le T,
\end{equation}
a.s. In other words, conditionally on the jump path, $\Delta^{Jb}_{A\to\widetilde A}$ belongs to $\cH_{\widetilde A,T}$.
\end{assumption}

\begin{remark}\label{rem:directionality}
The notation $A\to\widetilde A$ indicates that the discrepancy is represented in the Cameron--Martin space of the $\widetilde A$-Gaussian channel. If the Gaussian laws are equivalent, their Cameron--Martin spaces coincide as sets with equivalent norms, and a single representation such as \eqref{eq:CM_representation} is enough for equivalence of the full laws. If only one-sided Gaussian absolute continuity is known, the same representation gives the corresponding one-sided absolute continuity.
\end{remark}

\begin{theorem}[Gaussian noise with jumps]\label{thm:A-gen-main}
Let $A$ and $\widetilde A$ generate $C_0$-semigroups on $H$, and let $L=bt+W+Z$ with $W$ independent of $Z$. Suppose that
\begin{equation}\label{eq:Gaussian-one-sided-main}
        \Law(Y_A(W))\ll \Law(Y_{\widetilde A}(W))
        \quad\hbox{on }D_{H,T},
\end{equation}
and that Assumption~\ref{ass:CM} holds. Then
\[
        \Law(X^A)\ll \Law(X^{\widetilde A})
        \quad\hbox{on }D_{H,T}.
\]
If in \eqref{eq:Gaussian-one-sided-main} the two Gaussian laws are equivalent, then
\[
        \Law(X^A)\sim\Law(X^{\widetilde A})
        \quad\hbox{on }D_{H,T}.
\]
In particular, the inverse-covariance Hilbert--Schmidt condition \eqref{eq:Gaussian-HS}, together with Assumption~\ref{ass:CM}, implies equivalence of the two L\'evy-driven Ornstein--Uhlenbeck laws.
\end{theorem}

\begin{theorem}[Rigidity under pure jumps]\label{thm:B-gen-main}
Assume $L=Z$ is purely jump, with $W\equiv0$ and $b=0$. If
\[
        \Law(X^A)\ll\Law(X^{\widetilde A})
        \quad\hbox{on }D_{H,T},
\]
then $X^A=X^{\widetilde A}$ almost surely on $[0,T]$.
\end{theorem}

\section{Analytic semigroups and verifiable conditions}\label{sec:analytic}

We now specialise the abstract hypotheses to analytic semigroups. This is the setting most relevant to parabolic SPDEs generated by elliptic operators. The key point is that analytic smoothing converts fractional boundedness into the time-integrability required by the Gaussian comparison theorem.

\subsection{Fractional inverse-covariance boundedness implies Gaussian equivalence}

\begin{corollary}[Analytic semigroups]\label{cor:analytic_sufficient}
Let $A$ and $\widetilde A$ generate bounded analytic semigroups on $H$, and set $K:=\widetilde A-A$. Suppose that, for some $\beta\in(0,\tfrac12)$,
\begin{equation}\label{eq:R-beta-bounded}
        R_{A,\beta}:=Q^{-1/2}KA^{-\beta}
\end{equation}
is well-defined on $H$ and extends to an operator in $\cL(H)$. Then:
\begin{enumerate}[label=(\roman*)]
\item $KA^{-\beta}\in\cL(H)$, hence Duhamel's formula \eqref{eq:Duhamel} holds in the strong operator topology;
\item the inverse-covariance Hilbert--Schmidt condition \eqref{eq:Gaussian-HS} holds; consequently
\[
        \Law(Y_A(W))\sim\Law(Y_{\widetilde A}(W))
        \quad\hbox{on }D_{H,T}.
\]
\end{enumerate}
\end{corollary}

\begin{proof}
Since $Q^{1/2}$ is bounded, \eqref{eq:R-beta-bounded} implies
\[
        KA^{-\beta}=Q^{1/2}R_{A,\beta}\in\cL(H),
\]
which proves (i) by \Cref{thm:Duhamel}. For (ii), using \eqref{eq:analytic-smoothing},
\begin{align*}
\big\|Q^{-1/2}KS(t)Q^{1/2}\big\|_{\HS}
&=\big\|R_{A,\beta}A^\beta S(t)Q^{1/2}\big\|_{\HS}\\
&\le \|R_{A,\beta}\|\,\|A^\beta S(t)\|\,\|Q^{1/2}\|_{\HS}
 \le C t^{-\beta}.
\end{align*}
The square is integrable near $0$ precisely because $2\beta<1$.
\end{proof}

\begin{remark}\label{rem:resolvent-R}
In applications one may verify \eqref{eq:R-beta-bounded} through resolvent estimates of the form
\[
        \sup_{\lambda\in\Sigma_\theta}
        \big\|Q^{-1/2}K(\lambda+A)^{-1}\big\|\,|\lambda|^\beta<\infty,
\]
with the implicit requirement that the range of $K(\lambda+A)^{-1}$ is contained in $\Dom(Q^{-1/2})$. This is the covariance-weighted version of the usual sectorial functional-calculus criterion; see \cite[Ch.~II.4]{EngelNagel2000}.
\end{remark}

\subsection{Compound Poisson jumps}

For compound Poisson noise the Cameron--Martin condition can be checked pathwise. The statement below deliberately avoids an unconditional Novikov condition: after conditioning on the finitely many jumps, the Cameron--Martin shift is deterministic.

\begin{proposition}[Compound Poisson case]\label{prop:compound_Poisson}
Let $Z_t=\sum_{\tau_i\le t}\xi_i$ be a compound Poisson process with rate $\lambda>0$ and jump law $\nu$ on $H$, independent of $W$. Suppose that the hypotheses of \Cref{cor:analytic_sufficient} hold with some $\beta\in(0,\tfrac12)$. Then Assumption~\ref{ass:CM} holds on $[0,T]$. More explicitly,
\[
\Delta^{Jb}_{A\to\widetilde A}(t)
=\int_0^t\widetilde S(t-s)Q^{1/2}U^{A\to\widetilde A}_s\,ds,
\qquad 0\le t\le T,
\]
where $U^{A\to\widetilde A}$ is $\sigma(Z)$-measurable, independent of $W$, and belongs to $L^2(0,T;H)$ almost surely.
\end{proposition}

\begin{proof}
Put $R_{A,\beta}:=Q^{-1/2}KA^{-\beta}\in\cL(H)$. Since $Z$ is compound Poisson, the number $N(T)$ of jumps on $[0,T]$ is finite almost surely.

For the jump part,
\[
\Delta^J_{A\to\widetilde A}(t)
=\sum_{\tau_i\le t}\big(S(t-\tau_i)-\widetilde S(t-\tau_i)\big)\xi_i.
\]
By Duhamel's formula,
\[
S(r)-\widetilde S(r)=-\int_0^r\widetilde S(r-v)KS(v)\,dv,
\]
and therefore, after the change of variables $u=\tau_i+v$,
\begin{align*}
\Delta^J_{A\to\widetilde A}(t)
&=-\sum_{\tau_i\le t}\int_{\tau_i}^{t}
        \widetilde S(t-u)KS(u-\tau_i)\xi_i\,du\\
&=-\int_0^t\widetilde S(t-u)K
        \Big(\sum_{\tau_i<u}S(u-\tau_i)\xi_i\Big)\,du.
\end{align*}
The strict inequality $\tau_i<u$ is immaterial for the Lebesgue integral and avoids values at $u=\tau_i$. Define
\[
        F_J(u):=\sum_{\tau_i<u}A^\beta S(u-\tau_i)\xi_i,
        \qquad
        U^{(J)}_u:=-R_{A,\beta}F_J(u).
\]
Then $Q^{1/2}U^{(J)}_u=-KA^{-\beta}F_J(u)$ and hence
\begin{equation}\label{eq:jump-rep-revised}
        \Delta^J_{A\to\widetilde A}(t)
        =\int_0^t\widetilde S(t-u)Q^{1/2}U^{(J)}_u\,du.
\end{equation}

It remains only to check that $U^{(J)}\in L^2(0,T;H)$ a.s. For a fixed realisation with jumps $(\tau_i,\xi_i)_{i=1}^{N(T)}$, the analytic estimate \eqref{eq:analytic-smoothing} gives
\[
\int_{\tau_i}^{T}\|A^\beta S(u-\tau_i)\xi_i\|^2\,du
\le C\|\xi_i\|^2\int_0^T r^{-2\beta}\,dr<\infty.
\]
Since only finitely many jumps occur, $F_J\in L^2(0,T;H)$ a.s., and therefore $U^{(J)}\in L^2(0,T;H)$ a.s.

The deterministic drift part is treated similarly. Duhamel yields
\begin{align*}
\Delta^b_{A\to\widetilde A}(t)
&:=\int_0^t(S(t-s)-\widetilde S(t-s))b\,ds\\
&=-\int_0^t\widetilde S(t-u)K\Big(\int_0^uS(r)b\,dr\Big)\,du.
\end{align*}
Let
\[
        F_b(u):=A^\beta\int_0^uS(r)b\,dr,
        \qquad
        U^{(b)}_u:=-R_{A,\beta}F_b(u).
\]
Then
\begin{equation}\label{eq:drift-rep-revised}
        \Delta^b_{A\to\widetilde A}(t)
        =\int_0^t\widetilde S(t-u)Q^{1/2}U^{(b)}_u\,du.
\end{equation}
Moreover,
\[
        \|F_b(u)\|\le C\|b\|\int_0^u r^{-\beta}\,dr
        \le C'\|b\|u^{1-\beta},
\]
so $F_b\in L^2(0,T;H)$ and $U^{(b)}\in L^2(0,T;H)$.

Finally set
\[
        U^{A\to\widetilde A}:=U^{(J)}+U^{(b)}.
\]
Adding \eqref{eq:jump-rep-revised} and \eqref{eq:drift-rep-revised} gives \eqref{eq:CM_representation}. The process is a measurable functional of the compound Poisson path plus a deterministic term, hence is $\sigma(Z)$-measurable and independent of $W$.
\end{proof}

\begin{remark}\label{rem:compound-novikov}
An unconditional estimate of
\[
        \E\exp\left(\frac12\int_0^T\|U_s^{A\to\widetilde A}\|^2\,ds\right)
\]
is not part of the proposition and is not needed for \Cref{thm:A-gen-main}. A moment condition involving $\|Q^{1/2}\xi\|$ would control the wrong direction. The representative above involves $Q^{-1/2}KA^{-\beta}$ and the analytic smoothing of $S(\cdot)\xi$. Moreover, because the Poisson count is unbounded, a fixed exponential moment of $\|\xi\|^2$ alone does not by itself control the exponential of the square of the whole random shot-noise sum without further restrictions or localisation.
\end{remark}

\section{Counterexamples: failures of Cameron--Martin hypotheses}\label{sec:CMfail}

We record diagonal examples on $H=\ell^2(\N)$. They are intended to isolate possible failures of Cameron--Martin representability and of Novikov-type estimates. Except where explicitly stated, the examples are diagnostic and need not satisfy the Gaussian comparison hypothesis of \Cref{thm:GaussianEquiv}.

Let
\[
        S(t)e_n=e^{-a_nt}e_n,
        \qquad
        \widetilde S(t)e_n=e^{-\widetilde a_nt}e_n,
        \qquad
        Qe_n=q_ne_n,
\]
where $q_n>0$ and $\sum_nq_n<\infty$. Suppose that a single jump $\xi=\sum_n\xi_ne_n$ occurs at a deterministic time $s\in(0,T)$.

\begin{lemma}[Diagonal Cameron--Martin representative]\label{lem:diagonal_CM}
If \textup{(CM)}$_{A\to\widetilde A}$ holds for the single jump above, then the representative is
\[
        u_n(t)=\frac{\widetilde a_n-a_n}{\sqrt{q_n}}e^{-a_n(t-s)}\xi_n\,\mathbbm{1}_{[s,T]}(t),
\]
and
\[
        \|U\|_{L^2(0,T;H)}^2
        =\sum_{n=1}^\infty
        \frac{(\widetilde a_n-a_n)^2}{2a_nq_n}
        \big(1-e^{-2a_n(T-s)}\big)|\xi_n|^2.
\]
\end{lemma}

\begin{proof}
For $t\ge s$,
\[
\Delta^J_{A\to\widetilde A}(t)
=\sum_{n\ge1}\big(e^{-a_n(t-s)}-e^{-\widetilde a_n(t-s)}\big)\xi_ne_n.
\]
The $n$th coordinate of the Cameron--Martin identity is
\[
        \int_s^t e^{-\widetilde a_n(t-r)}\sqrt{q_n}\,u_n(r)\,dr
        =\big(e^{-a_n(t-s)}-e^{-\widetilde a_n(t-s)}\big)\xi_n.
\]
Differentiating after multiplying by $e^{\widetilde a_nt}$ gives the displayed expression for $u_n$, and integration in $t$ gives the norm formula.
\end{proof}

\begin{example}[Failure of the Cameron--Martin representative despite Gaussian smoothing]\label{ex:CM_fail_1}
Let $a_n=n^4$, $\widetilde a_n=n^4+1$, $q_n=e^{-n^2}$, and $\xi_n=e^{-n^2/4}/n$. Then
\[
        \sum_{n\ge1}\frac{(\widetilde a_n-a_n)^2}{a_nq_n}|\xi_n|^2
        =\sum_{n\ge1}\frac{e^{n^2/2}}{n^6}=\infty,
\]
so the Cameron--Martin representative is not in $L^2(0,T;H)$. On the other hand,
\[
        \sum_{n\ge1}\int_0^T(\widetilde a_n-a_n)^2e^{-2a_nt}\,dt
        \le\sum_{n\ge1}\frac1{2n^4}<\infty,
\]
which shows that Gaussian smoothing alone does not control arbitrary jump amplitudes in inverse-covariance norm.
\end{example}

\begin{example}[Asymmetric diagonal Cameron--Martin tests]\label{ex:CM_fail_2}
Let $a_n=n^4$, $\widetilde a_n\equiv1$, $q_n=n^{-8}$, and $\xi_n=n^{-7}$. For $A\to\widetilde A$ the terms in \Cref{lem:diagonal_CM} behave like
\[
        \frac{(\widetilde a_n-a_n)^2}{a_nq_n}|\xi_n|^2
        \asymp \frac1{n^2},
\]
so the representative is square-integrable. If the roles of $A$ and $\widetilde A$ are reversed, the corresponding terms behave like $n^2$, and the reversed representative fails to be in $L^2$. This example is outside the symmetric Gaussian comparison regime; it demonstrates that directional formulae must be checked in the reference Cameron--Martin channel actually used.
\end{example}

\begin{example}[Novikov may fail despite an $L^2$ representative]\label{ex:CM_fail_3}
In $H=\R$ with $Q=1$, let $a=\widetilde a+\delta$ for some $\delta>0$. If the jump size $\xi$ has a Student $t_\nu$ distribution with $\nu>2$, then $\E\xi^2<\infty$, so the Cameron--Martin representative has finite second moment by \Cref{lem:diagonal_CM}. However $\E e^{c\xi^2}=\infty$ for every $c>0$. Thus a Novikov-type exponential estimate is strictly stronger than $L^2$ representability and should not be built into a conditional Cameron--Martin criterion.
\end{example}

\begin{example}[Failure of a naive inverse-covariance factorisation]\label{ex:CM_fail_4}
Let $a_n\equiv1$, $\widetilde a_n=1+n^2$, and $q_n=n^{-6}$. A direct factorisation of $S-\widetilde S$ through the $\widetilde A$-Cameron--Martin channel would involve $Q^{-1/2}KS(t)$. But
\[
        \|Q^{-1/2}KS(t)\|_{\HS}^2
        =\sum_{n\ge1} n^6\,n^4e^{-2t}
        =e^{-2t}\sum_{n\ge1}n^{10}=\infty.
\]
Thus the inverse-covariance obstruction cannot be bypassed by a formal Duhamel factorisation.
\end{example}

\section{Proofs of the main theorems}\label{sec:proofs}

\subsection{Proof of Theorem~\ref{thm:A-gen-main}}
\begin{proof}
Let
\[
        \mu_A:=\Law(Y_A(W)),
        \qquad
        \mu_{\widetilde A}:=\Law(Y_{\widetilde A}(W)).
\]
For $z\in D_{H,T}$ define
\[
        c_A(z):=J_A(z)+B_A,
        \qquad
        c_{\widetilde A}(z):=J_{\widetilde A}(z)+B_{\widetilde A}.
\]
Then
\begin{equation}\label{eq:cA-cTilde}
        c_A(z)-c_{\widetilde A}(z)
        =\Delta^{Jb}_{A\to\widetilde A}(\cdot;z).
\end{equation}
By \Cref{prop:jump_reconstruction}, the maps $z\mapsto J_A(z)$ and $z\mapsto J_{\widetilde A}(z)$ are Borel on the full-measure set of jump paths; after redefining them arbitrarily outside this set, we regard them as Borel maps on $D_{H,T}$.

Disintegrating with respect to $Z$, the conditional law of $X^A$ given $Z=z$ is
\[
        \mu_A\circ\tau_{c_A(z)}^{-1},
\]
where $\tau_h(y)=y+h$. Similarly, the conditional law of $X^{\widetilde A}$ given $Z=z$ is
\[
        \mu_{\widetilde A}\circ\tau_{c_{\widetilde A}(z)}^{-1}.
\]
For $\Law(Z)$-a.e. $z$, Assumption~\ref{ass:CM} and \eqref{eq:cA-cTilde} give
\[
        c_A(z)-c_{\widetilde A}(z)\in\cH_{\widetilde A,T}.
\]
By \Cref{thm:CM_OU},
\begin{equation}\label{eq:shift-equiv-proof}
        \mu_{\widetilde A}\circ\tau_{c_A(z)}^{-1}
        \sim
        \mu_{\widetilde A}\circ\tau_{c_{\widetilde A}(z)}^{-1}
        \qquad\hbox{for }\Law(Z)\hbox{-a.e. }z.
\end{equation}
On the other hand, the Gaussian absolute continuity assumption \eqref{eq:Gaussian-one-sided-main} gives, for every fixed $z$,
\begin{equation}\label{eq:gaussian-shift-ac-proof}
        \mu_A\circ\tau_{c_A(z)}^{-1}
        \ll
        \mu_{\widetilde A}\circ\tau_{c_A(z)}^{-1},
\end{equation}
because translation by a fixed c\`adl\`ag path is a Borel isomorphism of $D_{H,T}$. Combining \eqref{eq:shift-equiv-proof} and \eqref{eq:gaussian-shift-ac-proof}, we obtain
\[
        \Law(X^A\mid Z=z)
        \ll
        \Law(X^{\widetilde A}\mid Z=z)
        \qquad\hbox{for }\Law(Z)\hbox{-a.e. }z.
\]
Integrating these conditional laws with respect to $\Law(Z)$ proves
\[
        \Law(X^A)\ll\Law(X^{\widetilde A}).
\]

If $\mu_A\sim\mu_{\widetilde A}$, then \eqref{eq:gaussian-shift-ac-proof} is an equivalence. Together with \eqref{eq:shift-equiv-proof} this yields equivalence of the conditional laws for $\Law(Z)$-a.e. $z$, and integration gives
\[
        \Law(X^A)\sim\Law(X^{\widetilde A}).
\]
\end{proof}

\subsection{Proof of Theorem~\ref{thm:B-gen-main}}
\begin{proof}
Assume $L=Z$. For a generator $B$, the jumps of $X^B=J_B(Z)$ are exactly the jumps of $Z$, because $S_B(0)=I$ and the stochastic convolution is continuous between jumps. Hence
\begin{equation}\label{eq:jump-reconstruct-solution}
        \mathcal Z(X^B)=Z
        \qquad\hbox{a.s.}
\end{equation}
Let
\[
        \Gamma_B:=\{x\in D_{H,T}:x=J_B(\mathcal Z(x))\}.
\]
By \Cref{prop:jump_reconstruction}, $\Gamma_B$ is Borel. Moreover, \eqref{eq:jump-reconstruct-solution} implies
\[
        \P(X^B\in\Gamma_B)=1.
\]
Since $\Law(X^A)\ll\Law(X^{\widetilde A})$ and $\Law(X^{\widetilde A})(\Gamma_{\widetilde A})=1$, we have
\[
        \P(X^A\in\Gamma_{\widetilde A})=1.
\]
But also $\P(X^A\in\Gamma_A)=1$. Therefore, using again $\mathcal Z(X^A)=Z$,
\[
        X^A=J_A(Z)=J_{\widetilde A}(Z)=X^{\widetilde A}
        \qquad\hbox{a.s. on }[0,T].
\]
\end{proof}

\section*{Further directions and open problems}\label{sec:open}

\begin{enumerate}[label=\textbf{P\arabic*.}, leftmargin=*, wide]

\item \textbf{Boundary case $\boldsymbol{\beta=\tfrac12}$.}
The analytic criterion uses $Q^{-1/2}KA^{-\beta}\in\cL(H)$ with $\beta<\tfrac12$, reflecting the integrability of $t^{-2\beta}$ near zero. It would be useful to identify additional spectral or covariance assumptions that allow the boundary case $\beta=\tfrac12$.

\item \textbf{Infinite-activity small jumps.}
For compound Poisson jumps the Cameron--Martin condition is pathwise because only finitely many jumps occur. For infinite-activity L\'evy measures one needs intrinsic conditions ensuring that the random series defining
\[
        Q^{-1/2}KA^{-\beta}\sum_{s<u}A^\beta S(u-s)\Delta Z_s
\]
belongs to $L^2(0,T;H)$ almost surely. Sharp near-zero and tail assumptions on the L\'evy measure remain to be identified.

\item \textbf{Degenerate Gaussian component.}
When $Q$ has non-trivial kernel, the relevant Cameron--Martin space is the range space of $Q^{1/2}$ and all drift perturbations must be projected onto this range. A complete comparison theorem should separate the controllable Gaussian directions from the jump-rigid directions.

\item \textbf{Banach-space extensions.}
It is natural to seek analogues in UMD Banach spaces, where Gaussian convolutions are described by $\gamma$-radonifying operators. The correct replacement for the inverse-covariance condition should be formulated in the corresponding reproducing-kernel Hilbert space.

\item \textbf{Non-analytic semigroups and $H^\infty$ calculus.}
Analytic smoothing was used only to obtain integrability of $A^\beta S(t)$ near zero. It would be interesting to replace analyticity by a bounded $H^\infty$ calculus, polynomial stability, or weighted square-function estimates.

\item \textbf{Time-dependent generators.}
For equations $dX(t)=A(t)X(t)\,dt+dL_t$ with an evolution family $U(t,s)$, one expects analogues of the Gaussian and Cameron--Martin criteria in terms of the two-parameter family $(t,s)\mapsto U(t,s)$.

\item \textbf{Quantitative Radon--Nikod\'ym bounds.}
Beyond qualitative equivalence, it would be valuable to estimate densities in terms of
\[
        \int_0^T\|Q^{-1/2}KS(t)Q^{1/2}\|_{\HS}^2\,dt
\]
and the size of the conditional Cameron--Martin representative of the jump--drift discrepancy.

\item \textbf{Numerical schemes and statistical identifiability.}
One may ask whether discretisations, such as Euler--Maruyama schemes with jumps, preserve approximate absolute continuity and how the inverse-covariance bounds affect statistical distinguishability of different generators.

\end{enumerate}

\section*{Acknowledgements}
The author is grateful to the referee for a careful reading of the manuscript and for pointing out the sign issue in the proof of the main comparison theorem and the inverse-covariance point in the compound Poisson argument. The author also thanks Grzegorz Bartosz for introducing him to the problem of comparing laws of L\'evy-driven Ornstein--Uhlenbeck processes.

\section*{Declarations}

\subsection*{Funding and/or Conflicts of Interest / Competing Interests}
The author declares that he has no known competing financial interests or personal relationships that could have appeared to influence the work reported in this paper.

\subsection*{Data Availability Statement}
No datasets were generated or analysed during the current study.


\begin{thebibliography}{99}

\bibitem{Applebaum}
D.~Applebaum,
\emph{L\'evy Processes and Stochastic Calculus}, 2nd ed.,
Cambridge Studies in Advanced Mathematics, vol.~116,
Cambridge Univ.\ Press, Cambridge, 2009.
\url{https://doi.org/10.1017/CBO9780511809781}

\bibitem{BartoszKania2019}
G.~Bartosz and T.~Kania,
Law equivalence of Ornstein--Uhlenbeck processes driven by a L\'evy process,
\emph{Indag.\ Math.\ (N.S.)} \textbf{30} (2019), no.~5, 796--804.
\url{https://doi.org/10.1016/j.indag.2019.05.003}

\bibitem{BrzVanN2000}
Z.~Brze\'zniak and J.M.A.M.~van Neerven,
Equivalence of Banach space-valued Ornstein--Uhlenbeck processes,
\emph{Stochastics Stochastics Rep.} \textbf{69} (2000), no.~1--2, 77--94.
\url{https://doi.org/10.1080/17442500008834233}

\bibitem{DaPratoZabczyk}
G.~Da Prato and J.~Zabczyk,
\emph{Stochastic Equations in Infinite Dimensions}, 2nd ed.,
Encyclopedia of Mathematics and its Applications, vol.~152,
Cambridge Univ.\ Press, Cambridge, 2014.
\url{https://doi.org/10.1017/CBO9781107295513}

\bibitem{EngelNagel2000}
K.-J.~Engel and R.~Nagel,
\emph{One-Parameter Semigroups for Linear Evolution Equations},
Graduate Texts in Mathematics, vol.~194,
Springer--Verlag, New York, 2000.
\url{https://doi.org/10.1007/b97696}

\bibitem{Kallenberg}
O.~Kallenberg,
\emph{Foundations of Modern Probability}, 2nd ed.,
Probability and its Applications,
Springer--Verlag, New York, 2002.
\url{https://doi.org/10.1007/978-1-4757-4015-8}

\bibitem{Loges1984}
W.~Loges,
Girsanov's theorem in Hilbert space and an application to the statistics of Hilbert space-valued stochastic differential equations,
\emph{Stochastic Process.\ Appl.} \textbf{17} (1984), no.~2, 243--263.
\url{https://doi.org/10.1016/0304-4149(84)90004-8}

\bibitem{Peszat1992}
S.~Peszat,
Equivalence of distributions of some Ornstein--Uhlenbeck processes taking values in a Hilbert space,
\emph{Probab.\ Math.\ Statist.} \textbf{13} (1992), no.~1, 7--17.

\bibitem{PeszatZabczyk2007}
S.~Peszat and J.~Zabczyk,
\emph{Stochastic Partial Differential Equations with L\'evy Noise: An Evolution Equation Approach},
Encyclopedia of Mathematics and its Applications, vol.~113,
Cambridge Univ.\ Press, Cambridge, 2007.
\url{https://doi.org/10.1017/CBO9780511721373}

\bibitem{Protter}
P.~E. Protter,
\emph{Stochastic Integration and Differential Equations}, 2nd ed.,
Springer, 2005.
\url{https://doi.org/10.1007/978-3-662-10061-5}

\end{thebibliography}
\end{document}